# A GENERALIZATION IN SPACE OF JUNG'S THEOREM

Florentin Smarandache
University of New Mexico
200 College Road
Gallup, NM 87301, USA
E-mail: smarand@unm.edu

In this short note we will prove a generalization of Joung's theorem in space.

**Theorem**. Let us have $n$ points in space such that the maximum distance between any two points is $a$. Prove that there exists a sphere of radius $r \leq a\frac{\sqrt{6}}{4}$ that contains in its interior or on its surface all these points.

*Proof:*
Let $P_1,...,P_n$ be the points. Let $S_1(O_1, r_1)$ be a sphere of center $O_1$ and radius $r_1$, which contains all these points. We note $r_2 = \max_{1 \leq i \leq n} P_i O_1 = P_1 O_1$ and construct the sphere $S_2(O_1, r_2)$, $r_2 \leq r_1$, with $P_1 \in Fr(S_2)$, where $Fr(S_2) =$ frontier (surface) of $S_2$.

We apply a homothety $H$ in space, of center $P_1$, such that the new sphere $H(S_2) = S_3(O_3, r_3)$ has the property: $Fr(S_3)$ contains another point, for example $P_2$, and of course $S_3$ contains all points $P_i$.

1) If $P_1, P_2$ are diametrically opposite in $S_3$ then $r_{\min} = \frac{a}{2}$.

If no, we do a rotation $R$ so that $R(S_3) = S_4(O_4, r_4)$ for which $\{P_3, P_2, P_1\} \subset Fr(S_4)$ and $S_4$ contains all points $P_i$.

2) If $\{P_1, P_2, P_3\}$ belong to a great circle of $S_4$ and they are not included in an open semicircle, then $r_{\min} \leq \frac{a}{\sqrt{3}}$ (**Jung**'s theorem).

If no, we consider the fascicle of spheres $S$ for which $\{P_1, P_2, P_3\} \subset Fr(S)$ and $S$ contains all points $P_i$. We choose a sphere $S_5$ such that $\{P_1, P_2, P_3, P_4\} \subset Fr(S_5)$.

3) If $\{P_1, P_2, P_3, P_4\}$ are not included in an open semisphere of $S_5$, then the tetrahedron $\{P_1, P_2, P_3, P_4\}$ can be included in a regulated tetrahedron of side $a$, whence we find that the radius of $S_5$ is $\leq a\frac{\sqrt{6}}{4}$.

If no, let's note $\max_{1\le i\le j\le 4} P_i P_j = P_1 P_4$. Does the sphere $S_6$ of diameter $P_1 P_4$ contain all points $P_i$?

If yes, stop (we are in the case 1).

If no, we consider the fascicle of spheres $S'$ such that $\{P_1, P_4\} \subset Fr(S')$ and $S'$ contains all the points $P_i$. We choose another sphere $S_7$, for which $P_5 \notin \{P_1, P_2, P_3, P_4\}$ and $P_5 \in Fr(S_7)$.

With these new notations (the points $P_1, P_4, P_5$ and the sphere $S_7$) we return to the case 2.

This algorithm is finite; therefore it constructs the required sphere.